\documentclass[a4paper,10pt]{article}

\usepackage{amssymb,amsmath,graphics}
\usepackage{bbm}
\newenvironment{dwd}{\par\noindent{\bf Proof.}}{\par\rightline{$\blacksquare$}}

\newtheorem{theo}{Theorem}
\newtheorem{prop}{Proposition}

\newtheorem{ex}[theo]{Example}
\newtheorem{lema}{Lemma}
\newtheorem{defi}{Definition}
\def\be#1\ee{\begin{equation}#1\end{equation}}
\newcommand{\ba}{\begin{eqnarray} }
\newcommand{\ea}{\end{eqnarray} }
\def\bt#1\et{\begin{theo}#1\end{theo}}
\def\bl#1\el{\begin{lema}#1\end{lema}}
\def\bp#1\ep{\begin{prop}#1\end{prop}}
\def\bd#1\ed{\begin{defi}#1\end{defi}}
\def\ccA{{\cal A}}
\def\ccB{{\cal B}}
\def\ccC{{\cal C}}

\def\ccM{{\cal M}}

\def\va{\varepsilon}
\def\n{\nu}

\def\E{\mathbf{E}}
\def\P{\mathbf{P}}

\def\R{{\mathbb R}}
\def\Z{{\mathbb Z}}

\def\ls{\leqslant}
\def\gs{\geqslant}

\def\1{{\mathbbm 1}}

\begin{document}

\title{\bf The suprema of infinitely divisible processes}
\author{{Witold Bednorz\footnote{Research supported by  NCN Grant UMO-2016/21/B/ST1/01489} and Rafa\l{} Martynek\footnote{Research supported by NCN Grant UMO-2018/31/N/ST1/03982.}}
\footnote{Institute of Mathematics, University of Warsaw, Banacha 2, 02-097 Warszawa, Poland}}
\date{}
\maketitle
\begin{abstract}
In this paper we complete the full characterization of the expected suprema of infinitely divisible processes. In particular, we remove the technical assumption called $H(C_{0},\delta)$ condition and settle positively the conjecture posed by M. Talagrand.
\end{abstract}

\section{Introduction}
\subsection{Infinitely divisible processes}
A stochastic process $X=(X_t)_{t\in T}$ over some abstract index set $T$ is called \textit{infinitely divisible} if all its finite dimensional marginal distributions are infinitely divisible (for the introduction to such distributions see e.g. \cite[Chapter 2]{Sat1}). An excellent reference for the study of general infinitely divisible processes containing all measure-theoretic subleties and examples is \cite{Ros2}. In particular, it contains some fundamental results concerning infinitely divisible processes, which we recap now.  First of all, an infinitely divisible process $X$ admits so-called \textit{L\'evy-Khintchine representation} of its characteristic function and as a consequence we can write that
$$X\overset{d}=G+Y,$$
where $G=(G_t)_{t\in T}$ is a centred Gaussian process and $Y=(Y_t)_{t\in T}$ is a \textit{Poissonian infinitely divisible process} independent of $G$ the definition of which we postpone a little. We refer to a Poissonian process as an infinitely divisible process without Gaussian component. Our main interest lies in the bounds of the expected supremum of $X$. It can be easily verified that the expected supremum of $X$ can be estimated from above and from below by the sum of expected suprema of $G$ and $Y$ using triangle inequality for the upper bound and Jensen's inequality for the lower bound. As we will recap in Section \ref{Overview} the suprema of Gaussian processes are very well-known, so it is sufficient to consider infinitely divisible processes without Gaussian component. Moreover, if $(X_t)_{t\in T}$ is infinitely divisible process and $(X_t')_{t\in T}$ is its' independent copy then $(X_t-X_t')_{t\in T}$ is also infinitely divisible. This procedure is called symmetrization and is of high importance for us since as we explain in Section \ref{A2} it allows us to assume with no loss in generality that the considered process is symmetric. Therefore, from now on we consider only symmetric infinitely divisible processes without Gaussian component. As a definition of such processes we treat a L\'evy-Khintchine representation of its characteristic function \cite[Theorem 2.8]{Ros2} which we present now.
\begin{defi}
Let $\R^{T}$ denote the space of functions from $T$ to $\R$. A stochastic process $(X_t)_{t\in T}$ is called \textit{real, symmetric, without Gaussian component infinitely divisible} if there exists a positive, $\sigma$-finite measure $\nu$ on $\R^T$ such that for each $t\in T$, $\int_{\R^T}(\beta(t)^2\wedge 1)\nu(d\beta)<\infty$, $\nu(\{0\})=0$ and for all families $(\alpha_t)_{t\in T}$ of real numbers and all $t\in T$ we have
\begin{equation}\label{ID}
\E\exp i\sum_{t\in T}\alpha_t X_t=\exp\left(-\int_{\R^T}\left(1-\cos\left(\sum_{t\in T}\alpha_t\beta(t)\right)\right)\nu(d\beta)\right).
\end{equation}
The measure $\nu$ is called the \textit{L\'evy measure} of the process.
\end{defi}

\noindent Let us emphasize an important note (\cite[p.334]{Tal1}) that the class of infinitely divisible processes is extremly large. In particular, the class of L\'evy processes on $\R^d$ is a subclass of infinitely divisible processes (for more examples see \cite{Ros2}). Moreover, let us notice that the fact that we are dealing with non-Gaussian infinitely divisible processes does not mean that they must be discontinuous. This is the case for L\'evy processes but there is a class of infinitely divisible processes like moving averages (see \cite{Ros3}) which are actually continuous.

\subsection{Rosi\'nski's representation and Poisson point process}
Usually in the theory of stochastic processes it is possible to discusss the sample boundedness (the question whether $\sup_{t\in T}X_t<\infty$ a.s.) of a process in terms of a single parameter - $\E\sup_{t\in T}X_t$ (see \cite{Tal4}). The example we have in mind are Gaussian processes for which the supremum in finite a.s. if and only if the expected supremum is finite. For many reasons this parameter is very convenient, though in general may not completely cover the question of sample boundedness. In particular, it is false for infinitely divisible processes (see Example \ref{ex} below). In this paper we prove the decomposition theorem (Theorem \ref{mainpop}) which completely characterize the role of the chaining method in the analysis of $\E\sup_{t\in T}X_t$. In this way, we explain the role of the chaining in the analysis of sample boundedness of infinitely divisible processes with
finite parameter $\E\sup_{t\in T}X_t$. 

\noindent From now on we assume that $\E\sup_{t\in T}X_t<\infty$. To achieve the goal of characterization of $\E\sup_{t\in T}X_t$ we use in the fundamental way the fact that infinitely divisible processes admit series representation. Let us define necessary objects to spell out the crucial result of Rosi\'nski (see \cite{Ros1}, \cite[Section 3.D]{Ros4}). 

\noindent Consider a $\sigma$-finite measure space $(\Omega,\nu)$ and a Poisson point process with intensity measure $\nu$ on $\Omega$ (see \cite{Dal} or \cite{Last} for the introduction to the subject of such processes). Recall it is a random subset $\Pi$ of $\Omega$ such that for any measurable subset $A$ of finite intensity measure the cardinality of $A\cap\Pi$
denoted by $|A\cap\Pi|$ is a Poisson random variable of expectation $\nu(A)<\infty$ and for disjoint measurable subsets $A_1,\dots, A_k$ the random variables $(|A_i\cap\Pi|)_{1\ls i\ls k}$ are independent. We denote the elements of $\Pi$ by $(Z_i)_{i\gs 1}$. Consider a Bernoulli sequence $(\va_i)_{i\gs 1}$ i.e. a sequence of independent random signs with $\P(\va_i=\pm1)=1/2$ independent of $Z_i$'s. It is easy to verify that a collection of random variables $(X_t)_{t\in T}$ given by
\begin{equation}\label{intr2}
X_t=\sum_{i\gs 1}\va_i Z_i(t).
\end{equation}
can be treated as infinitely divisible process in the sense of (\ref{ID}). Before we sketch the argument for that let us point out that it is more convenient to consider elements of $T$ as functions on $\R^{T}$ by identifying $t$ with a map $\beta\mapsto\beta(t)$ on $\R^{T}$ and write $t(Z_i)$ instead of $Z_i(t)$ in the above. 
\noindent We need the following elementary identity for integrable functions $f$ on $\Omega$ given by
\be\label{pois}
\E\sum_{i\gs1}f(Z_i)=\int_{\Omega}f(\omega)\nu(d\omega).
\ee
It follows instantly for indicators of sets and then for general $f$ by approximation. It is sometimes referred as Campbell's formula (\cite[Proposition 2.7]{Last}). Now, notice that by (\ref{pois}) we have that $\E\sum_{i\gs1}t(Z_i)^2\wedge1<\infty$ by the assumption that $\int_{\Omega}t^2\wedge 1 d\nu<\infty$. This applied conditionally on $Z_i$'s implies that $\sum_{i\gs 1}\va_i t(Z_i)$ converges almost surely.
Moreover, if $W$ is a Poisson random variable with mean $a$, then for $\lambda\in\R$
$$\E\exp(\lambda W)=\exp(a(\exp(\lambda)-1)).$$ 
Now, if $A$ is a measurable set of finite measure, $t=c\1_{A}$ and $W=|A\cap\Pi|$, then
$$\E\exp\left(\lambda\sum_{i\gs1}t(Z_i)\right)=\E\exp(\lambda cW)=\exp(\nu(A)(e^{\lambda c}-1)).$$
Therefore, by independence of $|A_i\cap\Pi|$ for step functions $t$ we have the following identity 
\be\label{expo}
\E\exp\left(\lambda\sum_{i\gs1}t(Z_i)\right)=\exp\left(\int(\exp(\lambda t(\omega))-1)\nu(d\omega)\right),
\ee
which extends to the bounded functions $t$ satisfying $\int |t|\wedge1d\nu<\infty$ by approximation. One can check using (\ref{expo}) that the characteristic function of linear combination of  $X_t$ given by (\ref{intr2}) satisfies (\ref{ID}) (see the computation at the beginning of Section \ref{Upper bound}). From now on we will refer to infinite divisible processes as to the collection $(X_t)_{t\in T}$ of random variables given by
\begin{equation}\label{intr3}
X_t=\sum_{i\gs 1}\va_i t(Z_i),
\end{equation}
where $T$ is a set of functions on $\Omega$ satisfying $\int_{\Omega}t^2\wedge 1 d\nu<\infty$ for $t\in T$. The following example returns to the problem of the boundedness and is using the series representation.

\begin{ex}\label{ex}
Let $\Omega=\R_+$, $\nu$ be the Lebesgue measure and $T=\{-t,t\}$ with $t(x)=x^{-2}$. Then $\sup_{s\in T}X_s=|X_t|<\infty$ a.s. since $\int t^2\wedge1 d\nu<\infty$, so $X_t$ is sample bounded. On the other hand by Khintchine's inequality
$$\E\sup_{s\in T}X_s\gs L\E\sqrt{\sum_{i\gs 1}t(Z_i)^2}\gs L\sup_{u\gs0}u\E\sqrt{\sum_{i\gs1}\1_{\{|t(Z_i)|\gs u\}}}\gs L\sup_{u\gs 1}u\nu(\{|t|\gs u\})=\infty,$$
where we used the fact that if $X$ is a Poisson random variable with mean $\lambda$, which is small then $\E\sqrt{X}$ is approximately $\lambda$ and $u\nu(\{|t|\gs u\})=\sqrt{u}\rightarrow\infty$ as $u\rightarrow\infty$.
\end{ex}

\subsection{Brief history of the problem}\label{History}
The main object of our study is $\E\sup_{t\in T}X_t$, where $X_t$ is given in (\ref{intr3}) and  
\begin{equation}\label{sup}
\E\sup_{t\in T}X_t = \sup\left\{\E\sup_{t\in F}X_t:\; F\subset T,\; F \;\mbox{finite}\right\}
\end{equation}
The question of the characterization of the expectation of the supremum of a stochastic process has a long history and goes back to the seminal results of Fernique and Talagrand on a centred Gaussian process $(G_t)_{t\in T}$ (in this section we abuse notation and  by $T$ we denote different index sets depending on the context) called the Majorising Measure Theorem (see \cite{Tal5} and \cite{Fer}).\\
\noindent Let us use notation $\lesssim$ for the inequality up to universal constant and $\asymp$, when both $\lesssim$ and $\gtrsim$ hold. We also write inequalities with constant $L$ which is a numerical constant but might be different at each occurance.\\ 
\noindent The classical formulation of the Majorising Measure Theorem states that 
$$\E\sup_{t\in T}G_t\asymp  \gamma_2'(T,d),$$
where for some distance $d$ on $T$ and the closed ball $B(t,\epsilon)$ centred at $t$ with radius $\epsilon>0$
$$\gamma_2'(T,d)=\inf_{\mu}\sup_{t\in T}\int_{0}^{\infty}\sqrt{\log\mu(B(t,\epsilon))^{-1}}d\epsilon.$$
The infimum is taken over all probability measures on $T$. The modern formulation of the Majorising Measure Theorem uses $\gamma_2$ functional which we define now. First, let $N_n=2^{2^n}$ for $n\gs1$ and $N_0=1$ and consider $T$ with some distance $d$. We will call nested sequence of partitions $(\ccA_n)_{n\gs0}$ of a set $T$ admissible if it holds that $|\ccA_0|=1$ together with $|\ccA_n|\ls N_n$ for $n\gs1$. By $A_n(t)$ we will denote (the unique) element of partition $\ccA_n$ that contains $t\in T$ and by $\Delta(\cdot)$ the diameter of a set in distance $d$. Define for $\alpha>0$
\begin{equation}\label{gam}
\gamma_{\alpha}(T,d)=\inf\sup_{t\in T}\sum_{n\gs0}2^{n/\alpha}\Delta(A_n(t)),
\end{equation}
where the infimum is taken over all admissible sequences. See \cite[Theorem 2.4.1]{Tal1} for the proof that $\E\sup_{t\in T}G_t\asymp  \gamma_2(T,d),$ where $d$ is the canonical distance $d(s,t)=(\E(X_t-X_s)^2)^{1/2}$. It can be actually verified directly that $\gamma_2'(T,d)\asymp\gamma_2(T,d)$ (see e.g. \cite{Tal3}). There is one more functional which is rather different from the other two. Namely,
$$\gamma''_2(T,d)=\sup_{\mu}\int_{T}\int_{0}^{\infty}\sqrt{\log\mu(B(t,\epsilon))^{-1}}d\epsilon\mu(dt),$$
where the supremum is over the probability measures on $T$. We will discuss those differences in the next section. What is important is that in the Gaussian case we have a single distance controlling the size of the increment i.e. for $u\gs 0$ we have that
$$\P(|G_t-G_s|\gs u)\ls 2\exp\left(-\frac{u^2}{d(s,t)^2}\right),$$
where $d(s,t)$ is the canonical distance given above. This situation is rather specific to the Gaussian case and it can be shown (see e.g. \cite{Bed2}) that in this setting
\begin{equation}\label{intuition}
\gamma_2(T,d)\asymp\gamma_2'(T,d)\asymp\gamma_2''(T,d).
\end{equation}
We know that a separable Gaussian processes admit a canonical representation called the Karhunen--Lo\`{e}ve representation (see e.g. \cite[Corollary 5.3.4]{M-R}). This observation was one of the motivations of the study of another canonical processes (see \cite[Chapter 10]{Tal1}) in particular the Bernoulli process defined as a collection $(B_t)_{t\in T}$, where now $T$ is a subset of $\ell^2$ and
\begin{equation}\label{ber}
B_t=\sum_{i\gs 1}\va_i t_i
\end{equation}
\noindent In the case of infinitely divisible processes there is a family of distances controlling the size of the increment, which makes the problem of characterization of $\E\sup_{t\in T}X_t$ fundamentally different. Initial attempt to deal with this task was done by M. Talagrand in 1993 \cite{Tal2}, where the regularity of infinitely divisible processes was proved under additional technical assumption on the L\'{e}vy measure. It was conjectured that the result holds without it, however little progress on this problem has been done since then. The main difficulty was related to the fact that the Bernoulli Conjecture (see \cite{Bed1}) was still open then. Nevertheless, many crucial properties of the Bernoulli process were proved in \cite{Tal2} like contraction principle or Sudakov's minoration principle, however the characterization of $\E\sup_{t\in T}B_t$ was missing. Under so-called $H(C_0,\delta)$ (see \cite[Definition 11.2.5]{Tal1}) condition on L\'evy measure $\nu$, which basically states that the measure $\nu$ is 'sufficiently anticoncentrated' it is shown in \cite{Tal2} that
\begin{equation}\label{charac}
\E\sup_{t\in T}X_t\asymp\inf\left\{\gamma_2(T_1,d_2)+\gamma_1(T_1, d_{\infty})+\E\sup_{t\in T_2}|X|_t: T\subset T_1+T_2 \right\},
\end{equation}
where $\gamma$-numbers are as in (\ref{gam}) for distances 
$$d_{\infty}(s,t)=\inf\{a>0: |s(\omega)-t(\omega)|\ls a\;\; \nu\mbox{-a.s.}\}$$
and 
$$d_2^2(s,t)=\int_{\Omega}(s(\omega)-t(\omega))^2\nu(d\omega)$$
and $|X|_t=\sum_{i\gs1}|t(Z_i)|$, where $T_2$ is a set functions $t$ on $\Omega$ satisfying $\int_{\Omega}|t|\wedge 1 d\nu<\infty$.\\

\noindent The characterization of $\E\sup_{t\in T}B_t$ given in \cite{Bed1} shed a new light on the result on $\E\sup_{t\in T}X_t$. Firstly, in one of the initial versions of \cite{Tal} the proof of (\ref{charac}) was significantly simplified by using to the fullest the result on the Bernoulli process. On the basis of this simplification, the proof of a complete characterization of $\E\sup_{t\in T}X_t$ without any additional assumptions on $\nu$ can be provided. The presentation of this result is the matter of the next section.

\section{Overiew of the results}\label{Overview}
The goal of this paper is to prove
\bt\label{mainpop}
Consider a countable set $T$ of measurable functions on $\Omega$ such that each $t\in T$ satisfies $\int_{\Omega}t^2\wedge 1 d\nu<\infty$ and assume that $0 \in T$. Let $X_t$ be as in (\ref{intr3}). Then 
$$\E\sup_{t\in T}X_t\asymp\inf\left\{\gamma_2(T_1,d_2)+\gamma_1(T_1, d_{\infty})+\E\sup_{t\in T_2}|X|_t: T\subset T_1+T_2,\;0\in T_1 \right\}.$$
\et
\noindent 
\noindent Observe that having the decomposition of $t\in T$ given by $t=t_1+t_2$, $t_1\in T_1, t_2\in T_2$ and $T\subset T_1+T_2$ we can write
$$\E\sup_{t\in T}X_t\ls\E\sup_{t\in T}X_{t_1}+\E\sup_{t\in T}X_{t_2}\ls\E\sup_{t\in T_1}X_t+\E\sup_{t\in T_2}X_t.$$
Theorem \ref{mainpop} states that an infinitely divisible process can be splitted into two parts. $\E\sup_{t\in T_1}X_t$ is the one explained through the chaining by ($\gamma_2(T_1,d_2)+\gamma_1(T_1, d_{\infty})$) and $\E\sup_{t\in T_2}X_t$ is bounded since $\E\sup_{t\in T_2}|X_t|\ls \E\sup_{t\in T_2}|X|_t$ and then by $\E\sup_{t\in T_2}|X|_t\ls L\E\sup_{t\in T}X_t.$
The two parts are of very different nature. The first part should be considered as the one where cancellations between terms occur while in the second there are no cancellations.

\noindent Before we continue, let us bring up a few classical facts used in the theory of suprema of stochastic processes to reassure that there is no loss generality in the assumptions of the above theorem. Firstly, $\E\sup_{t\in T}X_t$ is the same as taking supremum over set $T-s=\{t-s: t\in T\}$ for some $s\in T$, therefore we can always assume that $0\in T$. Note, that by the symmetry of $(X_t)_{t\in T}$ we can write
$$\sup_{s,t\in T}|X_s-X_t|=\sup_{s,t\in T}(X_s-X_t)=\sup_{s\in T}X_s+\sup_{t\in T}(-X_t),$$
so by taking expectatons we get that $\E\sup_{s,t\in T}|X_s-X_t|=2\E\sup X_t$. This combined with triangle inequality gives that for any $t_0\in T$ 
$$\E\sup_{t\in T}|X_t|\ls 2\E\sup_{t\in T}X_t+\E|X_{t_0}|\ls 3\E\sup_{t\in T}|X_t|.$$
In particular, when we assume that $0\in T$ then 
\be\label{zero}
\E\sup_{t\in T}|X_t|\asymp \E\sup_{t\in T}X_t.
\ee
We already mentioned that assuming the symmetry of $X_t$ is not any restriction in the context of suprema. Indeed, any process $(X_t)_{t\in T}$ can be replaced by the process $(X_t-X_t')_{t\in T}$, where $(X_t')_{t\in T}$ is an independent copy of $(X_t)_{t\in T}$. Section \ref{A2} discusses how to recover the result when $(X_t)$ is not symmetric, in Section \ref{A1} how the result can be extended to a separable set $T$. In Section \ref{A3} we provide an alternative formulation of Theorem \ref{mainpop} which might be more applicable. 
\subsection{Upper bound}\label{Upper bound}
An important feature of infinitely divisible processes is that they obey the Bernstein inequality. Denote by $\E_{\va}$ the expectation in $\va_i$'s. We have
$$\E_{\va}\exp\left(\lambda\sum_{i\gs1}\va_i t(Z_i)\right)=\exp\left(\sum_{i\gs 1}\log(\cosh(\lambda t(Z_i)))\right),$$
so by taking the expectation and applying (\ref{expo}) we get
$$\E\exp\left(\lambda\sum_{i\gs1}\va_i t(Z_i)\right)=\exp\left(\int(\cosh(\lambda t(\omega))-1)\nu(d\omega)\right).$$
Observe that $\cosh(\lambda t(\omega))-1\ls\lambda^2t^2(\omega)$ for $|\lambda t(\omega)|\ls 1$.
Consequently, arguing just as in the proof of \cite[Lemma 4.3.4]{Tal1}  we get
$$\P\left(\left|\sum_{i\gs1}\va_i(s(Z_i)-t(Z_i))\right|\gs v\right)\ls2\exp\left(-\frac{1}{L}\min\left(\frac{v^2}{d_2(s,t)^2},\frac{v}{d_{\infty}(s,t)}\right)\right).$$
\noindent This together with the chaining method \cite[ Theorem 2.2.23]{Tal1} yields

\be\label{Bern}
\E\sup_{t\in T}X_t\lesssim \gamma_2(T,d_2)+\gamma_1(T,d_\infty).
\ee
Obviously we also have that $|X_t|\ls\sum_{i\gs1}|t(Z_i)|$, so
\be\label{Gora}
\E\sup_{t\in T}X_t\ls \E\sup_{t\in T}|X|_t.
\ee
As discussed after the statement of Theorem \ref{mainpop} both (\ref{Bern}) and (\ref{Gora}) can be combined. The main result of this paper is that this combination can be reversed. We discuss the strategy for proving this in the next section.

\subsection{Strategy for the lower bound}
The starting point for obtaining the lower bound is the Gin\'e-Zinn type of result (\cite[Theorem 11.5.1]{Tal1}). We need to deal with the fact that $Z_i$'s are not independent, for which  the following lemma (cf. \cite[Proposition 3.8]{Last}) is used.
\begin{lema}\label{lemapoint}
Consider a Poisson point process of intensity $\nu$ and a measurable set $A$ with $0<\nu(A)<\infty$. Given $|A\cap\Pi|=N$, the set $A\cap\Pi$ has the same distribution as the set $\{X_1, X_2,\dots, X_N\}$, where the variables $X_i$ are independent and distributed according to the probability $\lambda$ on $A$ given by $\lambda(B)=\nu(A\cap B)/\nu(A)$ for a measurable subset $B$ of $A$. 
\end{lema}

\bt We have
\be\label{levy1}
\E\sup_{t\in T}|X|_t\ls\sup_{t\in T}\int_{\Omega}|t(\omega)|\nu(d\omega)+4\E\sup_{t\in T}|X_t|.
\ee
\et 
\begin{dwd}
Consider a measurable subset $A\subset\Omega$ with $\nu(A)<\infty$. Consider a sequence $\{X_1,\dots, X_N\}$ of independent random variables distributed according to the probability measure $\lambda$ on $A$ given for $B\subset A$ by $\lambda(B)=\nu(A\cap B)/\nu(A)$. As usual, by $(\va_i)_{i\gs1}$ we denote an independent Bernoulli sequence, independent of $X_i'$s. First,
$$\sum_{i\ls N}|t(X_i)|\ls\sum_{i\ls N}\E|t(X_i)|+\left|\sum_{i\ls N}(|t(X_i)|-\E|t(X_i)|)\right|$$
So,
\be\label{emp1}
\E\sup_{t\in T}\sum_{i\ls N}|t(X_i)|\ls\sup_{t\in T}\sum_{i\ls N}\E|t(X_i)|+\E\sup_{t\in T}\left|\sum_{i\ls N}(|t(X_i)|-\E|t(X_i)|)\right|
\ee
Observe that $\sum_{i\ls N}\E|t(X_i)|=N\nu(A)^{-1}\int_{A}|t(\omega)|\nu(d\omega)$. To deal with the second term consider an independent copy $(X_i')_{i\ls N}$ of $(X_i)_{i\ls N}$ and note that $(|t(X_i)|-|t(X_i')|)_{i\ls N}$ has the same law as $\va_i(|t(X_i)|-|t(X_i')|)_{i\ls N}$ by the independence of $X_i$'s. So, by Jensen's inequality
\begin{align*}
\E\sup_{t\in T}\left|\sum_{i\ls N}(|t(X_i)|-\E|t(X_i)|)\right|&\ls\E\sup_{t\in T}\left|\sum_{i\ls N}(|t(X_i)|-|t(X_i')|)\right|\\
&=\E\sup_{t\in T}\left|\sum_{i\ls N}\va_i(|t(X_i)|-|t(X_i')|)\right|.
\end{align*}
By triangle inequality
$$\E\sup_{t\in T}\left|\sum_{i\ls N}\va_i(|t(X_i)|-|t(X_i')|)\right|\ls2\E\sup_{t\in T}\left|\sum_{i\ls N}\va_i|t(X_i)|\right|.$$
The last step is to apply contraction principle (\cite[Theorem 5.3.6]{Tal1}) for a given value of $X_i$ to obtain
$$\E\sup_{t\in T}\left|\sum_{i\ls N}\va_i|t(X_i)|\right|\ls2\E\sup_{t\in T}\left|\sum_{i\ls N}\va_i t(X_i)\right|.$$
Hence we get
\be\label{GZinf2}
\E\sup_{t\in T}\sum_{i\ls N}|t(X_i)|\ls \frac{N}{\nu(A)}\sup_{t\in T}\int_{A} |t(\omega)|\nu(d\omega)+4\E\sup_{t\in T}\left|\sum_{i\ls N}\va_i t(X_i)\right|.
\ee
Now, we want to apply Lemma \ref{lemapoint}. For this consider a Poisson point process $(Z_i)_{i\gs1}$ of intensity measure $\nu$. Given $N=|\{i\gs1: Z_i\in A\}|$ we take the expectation in (\ref{GZinf2}) and use the fact that $\E N=\nu(A)$ to get
\be\label{GZinf3}
\E\sup_{t\in T}\sum_{i\ls N}|t(Z_i)|\1_A(Z_i)\ls \sup_{t\in T}\int_{A} |t(\omega)|\nu(d\omega)+4\E\sup_{t\in T}\left|\sum_{i\ls N}\va_i t(Z_i)\1_A(Z_i)\right|.
\ee
Obviously,
$$\int_{A} |t(\omega)|\nu(d\omega)\ls\int_{\Omega} |t(\omega)|\nu(d\omega).$$
For the second part we apply Jensen's inequality by taking the expectation in those $\va_i$'s at which $\1_A(Z_i)=0$ outside the absolute value and the supremum and get
$$\E\sup_{t\in T}\left|\sum_{i\ls N}\va_i t(X_i)\1_A(Z_i)\right|\ls\E\sup_{t\in T}\left|\sum_{i\ls N}\va_i t(X_i)\right|,$$
which finishes the proof by (\ref{GZinf3}). 
\end{dwd} 
\noindent The Gin\'e-Zinn inequality (\ref{levy1}) directs our attention to the essence of the problem. Notice that if there was no integral term on the right side of (\ref{levy1}) we could finish the proof. Therefore, we need a tool which provides a control over this term. This is given in the next theorem. To state it, let us define mappings which are square of distances, but we will refer to them as distances for simplicity. For $j\in\Z$, $r>0$ and $s,t\in T$ let
\be\label{dist}
\varphi_{j}(s,t)=\int_{\Omega}\left(r^{2j}|s(\omega)-t(\omega)|^2\wedge 1\right) \n(d\omega)
\ee
and we define for $t\in T$ and $\epsilon>0$,  $B_j(t,\epsilon)=\{s\in T:\varphi_j(s,t)\ls \epsilon\}$.
\bt\label{deco}\cite[Theorem 5.2.7]{Tal1}
Consider a countable set $T$ of measurable functions on $\Omega$, a number $r\gs4$ and assume $0\in T$. Consider an admissible sequence of partitions $(\ccA_n)_{n\gs0}$ of $T$, and for $A\in\ccA_n$ consider $j_n(A)\in\Z$, with the following properties, where $u>0$ is a parameter
$$A\in\ccA_n, B\in\ccA_{n-1}, A\subset B\implies j_n(A)\gs j_{n-1}(B),$$
$$\forall s,t\in A\in\ccA_n,\;\; \varphi_{j_n(A)}(s,t)\ls u2^n.$$
Then we can write $T\subset T'+T''+T'''$ where $0\in T'$ and
\be\label{gam2}\gamma_2 (T',d_2)\ls L\sqrt{u}\sup_{t\in T}\sum_{n\gs 0}2^n r^{-j_n(A_n(t))},
\ee
\be\label{gam1}
\gamma_1(T',d_{\infty})\ls L\sup_{t\in T}\sum_{n\gs 0}2^n r^{-j_n(A_n(t))}
\ee
\be\label{L1}
\forall t\in T'', \int_{\Omega}|t|d\nu \ls Lu\sup_{t\in T}\sum_{n\gs 0}2^n r^{-j_n(A_n(t))}.
\ee
Moreover,
\be\label{T3}
\forall t\in T''', \exists s\in T,\;\; |t|\ls 5|s|\1_{\{2|s|\gs r^{-j_0(A_0(t))}\}}.
\ee
\et
\noindent From now on, the parameter $r\gs4$ is fixed. What Theorem \ref{deco} effectively says is that if we can provide an admissible sequence of partitions of the set $T$ together with the sequence $(j_n)$ satisfying the described properties and most importantly that if 
\be\label{cel}
\sup_{t\in T}\sum_{n\gs 0}2^n r^{-j_n(A_n(t))}\ls\E\sup_{t\in T}X_t
\ee
then we are basically done with the proof of Theorem  \ref{mainpop}. Let us point out that providing an admissible sequence of partitions such that the series $\sum_{n\gs0}2^n r^{-j_n(A_n(t))}$ is a lower bound for some functional related to $\E\sup_{t\in T}X_t$ is a central tool in Talagrand's machinery (see \cite[Theorem 10.1.2]{Tal1}). However, the application of this theorem is only possible with $H(C_0,\delta)$ condition additionally assumed. Without it it seems hopeless. Therefore, we have to engage another strategy which we outline now.\\

\noindent The fundamental idea behind the argument is that conditionally on $(Z_i)_{i\gs 1}$, $(X_t)_{t\in T}$ is a Bernoulli process. Then, the idea is to formulate analogues of functionals $\gamma_2$, $\gamma_2'$ and $\gamma_2''$ introduced in Section \ref{History} in the Bernoulli case. The importance of functionals $\gamma_2$ and $\gamma_2'$ are of different nature than $\gamma_2''$. Notice, that if we consider a random distance $d=d_{\omega}$ (consider a conditionally Gaussian process), then it can be shown that if $d_{\omega}$ concentrates well around its mean it follows that $\E\gamma_2''(t, d_\omega)\gtrsim\gamma_2''(T, \E d_{\omega})$. This type of property cannot be obtained with $\gamma_2$ or $\gamma_2'$ because the infimum which is present in the definition of both of these functionals cannot be pulled out of the expectation. On the other hand, the functionals $\gamma_2$ or $\gamma_2'$ as opposed to $\gamma_2''$ are natural to the chaining arguments. What lies in the heart of the proof of Theorem \ref{mainpop} is finding  the "correct" analogues of $\gamma_2$, $\gamma_2'$ and $\gamma_2''$ in the setting of Bernoulli processes and proving a version of (\ref{intuition}) namely that 
$$\gamma_2\lesssim\gamma_2'\lesssim\gamma_2''$$
in the Bernoulli case. Now, we proceed to define the analogues of chaining functionals in the Bernoulli case. Unsurprisingly, they are much more complicated than in the Gaussian case since there are many distances involved. 

\noindent The random distance we consider is given by$$\tilde{\varphi}_{j,Z}(s,t):=\sum_{i\gs1}(r^{2j}|t(Z_i)-s(Z_i)|^2)\wedge 1.$$
Notice that by (\ref{pois}) we have that $\E\tilde{\varphi}_{j,Z}(s,t)=\varphi_j(s,t)$, where $\varphi_j(s,t)$ is as in (\ref{dist}).

\noindent For the moment let us focus on the purely Bernoulli setting, i.e. consider a subset $S$ of $\ell^2$ and the sequence of distances
$$\tilde{\varphi}_{j}(s,t):=\sum_{i\gs1}(r^{2j}|t_i-s_i|^2)\wedge 1.$$
Let $k_0$ be the largest integer such that the diameter of $S$ in the $\ell^2$-distance does not exceed $r^{-k_0}$. For $t\in S$ define $k_0(t)=k_0$ and for $n\gs 1$ and fixed probability measure $\mu$ on $S$
\be\label{k}
k_n(t)=\sup\{j\in\Z\;;\;\mu(\{s\in T, \tilde{\varphi}_j(s,t)\ls 2^n\})\gs N_{n}^{-1}\}.
\ee
Let 
\be\label{func1}
I_{\mu}(t)=\sum_{n\gs0}2^n r^{-k_n(t)}.
\ee 
Define 
$$\beta'(S,(\tilde{\varphi_j})):=\inf_{\mu}\sup_{t\in S}I_{\mu}(t),\;\;\;\;\;\;\;\;\;\;\       \beta''(S, (\tilde{\varphi_j})):=\sup_{\mu}\int_{S}I_{\mu}(t)\mu(dt).$$
The critical result of Talagrand \cite{Tal} states that 
\be\label{critical}
\beta''(S, (\tilde{\varphi_j}))\lesssim\E\sup_{t\in S}B_t.
\ee
It gives a rise to the roadmap of inequalities leading to the proof of Theorem \ref{mainpop}. Before we present it let us define 
$$\beta(T, (\varphi_j)):=\inf\sup_{t\in T}\sum_{n\gs 0}2^n r^{-j_n(A_n(t))},$$ where the infimum runs over admissible sequences of partitions and labels $j_n$ from Theorem \ref{deco}. We have to argue that the following steps hold true.
$$\beta(T,(\varphi_{j}))\overset{(4)}{\lesssim}\beta'(T,(\varphi_j))\overset{(3)}{\lesssim}\beta''(T,(\varphi_j))\overset{(2)}{\lesssim}\E\beta''(T,(\tilde{\varphi}_{j,Z}))\overset{(1)}{\lesssim}\E\sup_{t\in T}X_t$$
and
$$\inf\left\{\gamma_2(T_1,d_2)+\gamma_1(T_1, d_{\infty})+\E\sup_{t\in T_2}|X|_t: T\subset T_1+T_2,\;0\in T_1 \right\}\overset{(5)}{\lesssim}\beta(T, (\varphi_j))+\E\sup_{t\in T}X_t.$$
In the above we understand $\beta'$ and $\beta''$ defined for the set $T$ are almost just as in the Bernoulli case, but we postpone formal definitions to keep the above picture intuitive. Steps $(1)$, $(2)$ were communicated to authors by M. Talagrand as one of the preprints of \cite{Tal}, step $(5)$ is essentialy present in \cite{Tal1}. Steps $(3)$ and $(4)$ are the main contribution of this paper. We provide details of each step in the following sections to make this paper self-contained. Before that, let's summarize each of them. Step $(1)$ is the major application of the characterization of Bernoulli processes \cite{Bed1}, step $(2)$ is the analogue of the discussed procedure of pulling the supremum out of expectation  when dealing with random  distance ($\E\gamma_2''(t, d_\omega)\gtrsim\gamma_2''(T, \E d_{\omega})$), step $(3)$ is a consequence of minimax theorem, step $(4)$ requires building the admissible sequence of partitions together with labels $j_n$ satisfying the conditions of Theorem \ref{deco}. Finally, step $(5)$ might be already clear as a consequence of Gin\'e-Zinn inequality (\ref{GZinf3}) and Theorem \ref{deco} with a minor issue of the term in (\ref{T3}) which we have to take care of.

\section{Proofs}

\subsection{Proof of Inequality (1)}

We have the following result due to the first named author and R. Lata\l{}a called the Bernoulli Theorem \cite{Bed1}.

\bt\label{Ber}
Let $S\subset\ell^2$. Then it holds that
$$\inf\{\gamma_2(S_1,d_{\ell^2})+\sup_{t\in S_2}\|t\|_1;S\subset S_1+S_2\}\asymp \E\sup_{t\in S}B_t,$$
where $d_{\ell^2}$ denotes the $\ell^2$-distance and $\|\cdot\|_1$ is the $\ell^1$-norm.
\et 
\noindent 
Our main goal in this section is to prove $(\ref{critical})$ from which Inequality $(1)$ follows easily. Define 
$$b^{*}(S):=\inf\{\gamma_2(S_1,d_{\ell^2})+\sup_{t\in S_2}\|t\|_1;S\subset S_1+S_2\},$$ 
$$\beta(S,(\tilde{\varphi}_j)):=\inf\sup_{t\in S}\sum_{n\gs 0}2^n r^{-j_n(A_n(t))},$$ where the infimum is taken over admissible sequences $(\ccA_n)_{n\gs0}$ of partititons of $S$ and  integers $j_n(A)$ for $A\in\ccA_n$ satisfying
\be\label{func2}
s,t\in A\implies \tilde{\varphi}_{j_n(A)}(s,t)\ls 2^n
\ee
and
$$\Delta(S,d_{\ell^2})\ls r^{-j_0(S)}.$$

\noindent The proof of (\ref{critical}) relies on the two following inequalities
$$\sup_{\mu}\int_S I_{\mu}(t)\mu(dt)\overset{(B2)}{\lesssim}\beta(S,(\tilde{\varphi}_j))\overset{(B1)}{\lesssim}b^{*}(S)$$  
combined with Theorem \ref{Ber}.
Inequality $(B1)$ can be found in \cite[Proposition 5.2.5]{Tal1}. The novelty of \cite{Tal} is the Inequality $(B2)$ which we include together with the proof.
\begin{lema}
Given any probability measure $\mu$ on $S$ we have
\be\label{step2}
\int_S I_{\mu}(t)\mu(dt)\ls L\beta(S,(\tilde{\varphi}_j)).
\ee
\end{lema}
\begin{dwd}
Consider an admissible sequence $(\ccA_n)_{n\gs0}$ of partitions of $S$ and for $A\in\ccA_n$ and integers $j_n(A)$ as in (\ref{func2}) such that 
$$\sup_{t\in S}\sum_{n\gs0}2^n r^{-j_n(A_n(t))}\ls2\beta(S,(\tilde{\varphi}_j)).$$
By the definition of $k_0$ it follows that $\tilde{\varphi}_{k_0+1}(s,t)>1$ for some $s,t\in S$ and since $\tilde{\varphi}_{k_0+1}(s,t)\ls r^{2(k_0+1)}d_{\ell^2}(s,t)^2$ we have $r^{-k_0-1}<\Delta(S,d_{\ell^2})$. Hence, since $\Delta(S,d_{\ell^2})\ls r^{-j_0(S)}$,
$$r^{-k_0-1}\ls r^{-j_0(S)}\ls \sup_{t\in S}\sum_{n\gs0}2^n r^{-j_n(A_n(t))}.$$
Now, for $n\gs1$, $A\in\ccA_n$ and $t\in A$
$$A\subset\{s\in S:\;\tilde{\varphi}_{j_n(A)}(s,t)\ls2^n\}\subset\{s\in S:\;\tilde{\varphi}_{j_n(A)}(s,t)\ls2^{n+1}\}.$$
So, if $\mu(A)\gs N_{n+1}^{-1}$ then $k_{n+1}(t)\gs j_{n}(A)$, hence
$$\int_A 2^{n+1}r^{-k_{n+1}(t)}\mu(dt)\ls2\int_A 2^n r^{-j_n(A_n(t))}\mu(dt).$$
On the other hand, if $\mu(A)<N_{n+1}^{-1}$
$$\int_A 2^{n+1}r^{-k_{n+1}(t)}\mu(dt)<2^{n+1}r^{-k_0}N_{n+1}^{-1}.$$
Summation over $A\in\ccA_n$ and then over $n\gs 0$ yields
$$\int_S\sum_{n\gs1}2^n r^{-k_n(t)}\ls L\sup_{t\in S}\sum_{n\gs0}2^n r^{-j_n(A_n(t))}+Lr^{-k_0},$$
which finishes the proof, because for $n=0$, $$2^n r^{-k_n(t)}\ls Lr^{-k_0}\ls L\sup_{t\in S}\sum_{n\gs0}2^n r^{-j_n(A_n(t))}.$$
\end{dwd}
\noindent Having prepared the lower bound for Bernoulli processes, we are ready to come back to infinitely divisible processes. Consider $I_{\mu, Z}$ defined as in (\ref{func1}) but for $s, t\in T$, $$\tilde{\varphi}_{j, Z}(s,t)=\sum_{i\gs1}(r^{2j}|t(Z_i)-s(Z_i)|^2)\wedge 1$$ and random numbers $(k^{Z}_{n})_{n\gs 1}$ defined as in (\ref{k}) but with $\tilde{\varphi}_{j, Z}(s,t)$
and
$$k_0^Z=\sup\{j\in\Z:\;\forall s,t\in T,\;\tilde{\varphi}_{j, Z}(s,t)\ls 1\}\in\Z\cup\{\infty\}.$$ 
Conditionally on $Z_i$'s, $X_t$ is a Bernoulli process. Thus, if we denote by $\E_{\va}$ the expectation with respect to $\va_i$'s, then by (\ref{critical}) we get for any probability measure $\mu$ on $T$ that
$$\int_{T}I_{\mu, Z}(t)\mu(dt)\ls L\E_{\va}\sup_{t\in T}X_t$$ 
By taking expectation with respect to  $Z_i$'s  we obtain Inequality (1) 
\be\label{ID1}
\E\sup_{\mu}\int_{T}I_{\mu,Z}(t)\mu(dt)\ls L\E\sup_{t\in T}X_t.
\ee
\subsection{Proof of Inequality (2)}
Using the distance (\ref{dist}) we define numbers $(j_k^{\mu}(t))_{k\gs0}$ associated to each point $t\in T$ and probability measure $\mu$ on $T$. Let
\be\label{j0}
j_0=\sup\{j\in\Z: \forall s,t\in T, \varphi_j(s,t)\ls 4\}.
\ee
The use of constant $4$ in the above will become apparent with the next two results.
Given any probability measure $\mu$ on $T$ we define for any integer $n\gs0$ and $t\in T$
$$
j_0^\mu(t)=j_0
$$
and for $n\gs 1$
\be\label{jk}
j_n^{\mu}(t)=\sup\{j\in\Z: \mu(B_{j}(t,2^n))\gs N_{n}^{-1}\}.
\ee
We also define 
$$J_{\mu}(t)=\sum_{n\gs 0}2^n r^{-j_n^{\mu}(t)}.$$

\begin{lema}
For every $s,t\in T$ and any $j\in\Z$ the following inequality holds true.
\be\label{exp}
\P(\tilde{\varphi}_{j, Z}(s,t)\ls\varphi_{j}(s,t)/4)\ls\exp(-\varphi_j(s,t)/4).
\ee
\end{lema}

\begin{dwd}
Let $f(Z_i)=r^{2j}|Z_i(s)-Z_i(t)|^2\wedge 1$. Recall that by (\ref{pois}), $\E\tilde{\varphi}_{j,Z}(s,t)=\varphi_j(s,t)$. Notice that for $0\ls x\ls 1$ we have $e^{-x}\ls 1-x/2$, so due to (\ref{expo}) we have
\begin{align*}
\E\exp\left(-\sum_{i\gs1}f(Z_i)\right)&=\exp\left(\int(\exp(-f(\omega))-1)\nu(d\omega)\right)\\
&\ls\exp\left(-1/2\int f(\omega)\nu(d\omega)\right).
\end{align*}
For any constant $A$ we have
$$\P(\sum_{i\gs1}f(Z_i)\ls A)\ls \exp(A)\E\exp(-\sum_{i\gs1}f(Z_i)),$$ so the result follows by taking $A=\varphi_j(s,t)/4$
\end{dwd}
\noindent The following is the proof of Inequality $(2)$. Let us emphasize the fact that numbers $k_n^{Z}$ are random objects as opposed to $j_n^{\mu}$.
\begin{theo}\label{joty}
For each $t\in T$ we have $J_{\mu}(t)\ls L\E I_{\mu, Z}(t)$.
\end{theo}
\begin{dwd}
Fix a probability measure $\mu$ on $T$. To simplify the notation put $j_n(t)=j_n^{\mu}(t)$. We aim to prove that 
\be\label{exp1}
\P(k_0^Z\ls j_0)\gs\frac{1}{2}
\ee
and for $n\gs3$
\be\label{exp2}
\P(k_{n-2}^Z(t)\ls j_n(t))\gs\frac{1}{2}.
\ee
These relations imply respectively that $\E r^{-k_0^Z}\gs r^{-j_0}/L$ and for $n\gs3$, $\E2^{n-3}r^{-k_{n-3}^Z(t)}\gs2^{n}r^{-j_n(t)}/L$.
Summing these inequalities yields the result.\\
\noindent For (\ref{exp1}) notice that if $j_0=\infty$ there is nothing to prove. Otherwise, we use (\ref{exp}). By the definition of $j_0$, there exist $s,t\in T$ such that $\varphi_{j_0+1}(s,t)>4$ and therefore
$$\P(\tilde{\varphi}_{j_0+1,Z}(s,t)>1)\gs1-\exp(-1).$$
By the definition of $k_0^Z$ the event $\tilde{\varphi}_{j_0+1,Z}(s,t)>1$ implies that $k_0^Z\ls j_0$.\\
Now we argue that (\ref{exp2}) holds. If $j_n=\infty$ there is nothing to prove. By the definition of $k_n^Z(t)$ we have
$$\mu(\{s\in T:\; \varphi_{k_n^Z(t)+1, Z}(s,t)\ls 2^n\})\ls N_{n}^{-1}.$$
On the other hand, by (\ref{exp}), $\varphi_{k_n^Z(t)+1}(s,t)\gs 2^n$ implies that
$$\P(\tilde{\varphi}_{k_n^Z(t)+1, Z}(s,t)\ls 2^{n-2})\ls\exp(-2^{n-2})\ls N_{n-2}^{-1}.$$
Hence, 
$$\E\mu(\{s\in T:\;\varphi_{k_n^Z(t)+1}(s,t)\gs 2^n, \tilde{\varphi}_{k_n^Z(t)+1, Z}(s,t)\ls 2^{n-2}\})\ls N_{n-2}^{-1}$$
and by Markov inequality with probability $\gs 1/2$ we have
$$\mu(\{s\in T:\;\varphi_{k_n^Z(t)+1}(s,t)\gs 2^n, \tilde{\varphi}_{k_n^Z(t)+1, Z}(s,t)\ls 2^{n-2}\})\ls 2N_{n-2}^{-1}.$$
Therefore, conditioned on this event we get that 
$$\mu(\{s\in T:\;\tilde{\varphi}_{k_n^Z(t)+1, Z}(s,t)\ls 2^{n-2}\})\ls N_{n}^{-1}+2N_{n-2}^{-1}<N_{n-3}^{-1}$$
and in turn $k_{n-3}^Z(t)\ls j_n(t)$. This finishes the proof.
\end{dwd}
\noindent Theorem \ref{joty} together with (\ref{ID1}) imply that for any probability measure $\mu$
\be\label{ID2}
\int_{T}J_{\mu}(t)\mu(dt)\ls L\int_{T}\E I_{\mu, Z}(t)\mu(dt)=L\E\int_{T} I_{\mu, Z}(t)\mu(dt)\ls L\E\sup_{t\in T}X_t.
\ee
\subsection{Proof of Inequality (3)}
The inequality $\beta'(T,(\varphi_j))\lesssim\beta''(T,(\varphi_j))$ is in spirit the convexity argument already present in \cite[Proposition 3.2]{Tal4}. For this, we have to work on any finite subset of $T$ which we fix now and denote by $F$. Note that we do not prove Inequality $(3)$ for $T$ but for $F$, but this will be enough to proceed with Inequality $(4)$. We want to prove that there exists some probability measure $\mu_0$ on $F$ such that for each $t\in F$ we have
$$J_{\mu_0}(t)\ls\sup_{\mu}\int_{F}J_{\mu}(s)\mu(ds) $$
This could be done by a standard argument using Hahn-Banach theorem provided that $J_{\mu}(t)$ is convex as a function of $\mu$. Proving the convexity is a highly non-trivial task which is therefore replaced by a modified reasoning.  What we show is that $J_{\mu}$ is "essentially convex" which we describe in the next result.

\begin{theo}\label{essconv}
Consider any non-negative numbers $(\alpha_i)_{i\gs1}$ with $\sum_{i\gs 1}\alpha_i=1$ and any probability measures $(\mu_i)_{i\gs 1}$ on $F$. Then for any $t\in F$
$$J_{\sum_{i\gs 1}\alpha_i\mu_i}(t)\ls L\sum_{i\gs 1}\alpha_i J_{\mu_i}(t).$$ 
\end{theo}
\begin{dwd}
Fix $t\in F$.
Define numbers $j_n(t)$ with $j_0$ as in (\ref{j0}) and for $n\gs1$
\be\label{jnpop}
r^{-j_n(t)-1}<\sum_{i\gs 1}\alpha_i r^{-j_n^{\mu_i}(t)}\ls r^{-j_n(t)}.
\ee
It is a matter of changing the order of summation to notice that  
\be\label{mainlow}
\sum_{n\gs0}2^n r^{-j_n(t)}\ls L\sum_{i\gs 1}\alpha_i\sum_{n\gs0}2^nr^{-j_n^{\mu_i}(t)}=L\sum_{i\gs 1}\alpha_i J_{\mu_i}(t).
\ee
since $j_0$ is fixed. What is crucial is that $j_n$'s defined in (\ref{jnpop}) preserve the defining property given by (\ref{jk}) and indeed the quantity on the left hand side of inequality (\ref{mainlow}) is bounded below by $J_{\mu}(t)$ with $\mu=\sum_{i\gs 1}\alpha_i\mu_i$. Namely, we show that for $n\gs1$
\begin{equation}\label{jnpopmu}
\mu(B_{j_n(t)}(t,2^n))\gs\frac{2}{3}\frac{1}{N_n}>\frac{1}{N_{n+1}},
\end{equation}
which implies that $j_{n+1}^{\mu}(t)\gs j_n(t)$.\\
\noindent Observe that if $j^{\mu_i}_n (t)\gs j_n(t)$, then $\mu_{i}(B_{j_n(t)}(t,2^n))\gs1/N_n$. Define $\beta_j=\sum_{i}\alpha_i \1_{\{j_n^{\mu_i}(t)=j\}}$. Certainly, $\sum_{j\in\Z}\beta_j=1$. Moreover,
$$\mu(B_{j_n(t)}(t,2^n))\gs\sum_{j\gs j_n(t)}\beta_j\frac{1}{N_n}.$$
Now, we will argue that $\sum_{j<j_n(t)}\beta_j$ can be bounded from above so that (\ref{jnpopmu}) follows. From the definition (\ref{jnpop}) of $j_n(t)$ we have
$$\sum_{j\in\Z}\beta_j r^{-j}\ls r^{-j_n(t)},$$
which implies that $\beta_j r^{-j}\ls r^{-j_n(t)}$ and in turn that $\beta_j\ls r^{-j_n(t)+j}$. Hence, since $r\gs4$,
$$\sum_{j<j_n(t)}\beta_j\ls\sum_{j<j_n(t)}r^{j-j_n(t)}\ls\sum_{l\gs1}r^{-l}\ls\frac{1}{3},$$
which finishes the argument.
\end{dwd}

\bt\label{sep}
There exists a probability measure $\mu_0$ on $F$ such that for each $t\in F$ and we have
$$J_{\mu_0}(t)\ls L\sup_{\mu}\int_{F}J_{\mu}(s)\mu(ds)$$
\et
\begin{dwd}
It is almost a straightforward consequence of the minimax theorem (see e.g. \cite{von}) used for sets from which one is compact and the other is convex and then combined with Theorem (\ref{essconv}). We only need to be a little careful to make sure that we consider $J_{\mu}(t)$ which are finite for every $t\in F$. For this, let $\ccM$ be the class of probability measures $\mu$ on $F$ satisfying $$\mu(\{t\})\gs\frac{1}{2|F|}.$$
Notice that for $\mu\in\ccM$, $J_{\mu}(t)<\infty$. Let
$$\ccC=\mbox{conv}\left\{(J_{\mu}(t))_{t\in F}\in\R^{F}:\mu\in\ccM\right\}.$$
For $f\in\ccC$, $\sup_{\mu\in\ccM}\int_{F}fd\mu\gs\sup_{t\in F}f(t)/2$.
So, by the minimax theorem
$$\inf_{f\in\ccC}\sup_{t\in F}f(t)\ls2\inf_{f\in\ccC}\sup_{\mu\in\ccM}\int_{F}fd\mu=2\sup_{\mu\in\ccM}\inf_{f\in\ccC}\int_{F}fd\mu\ls2\sup_{\mu\in\ccM}\int J_{\mu}d\mu.$$
Hence, there is an element of $\ccC$ bounded from above by $2\sup\int J_{\mu}d\mu$, where now the supremum is over all probability measures on $F$. Namely, we get non-negative numbers $(\alpha_i)_{i\gs1}$ summing to $1$ and probability measures $\mu_i\in\ccM$ such that 
$$\sum_{i\gs 1}\alpha_i J_{\mu_i}(t)\ls 2\sup_{\mu}\int J_{\mu}d\mu.$$
The result follows from Theorem \ref{essconv} with $\mu_0=\sum_{i\gs 1}\alpha_i\mu_i$.

\end{dwd}

\noindent 
\subsection{Proof of Inequality (4)}
What we achieved in the previous section is the lower bound from a majorising measure $\mu$ which depends on the finite set $F$. Namely, 
$$\sup_{t\in F}J_{\mu}(t)\ls L\E\sup_{t\in T}X_t.$$
We aim to replace the lower bound by the sum depending on the admissible sequence of partitions of the whole set $T$ as in Theorem \ref{deco}, which is precisely the meaning of Inequality $(4)$.
Before stating the result observe (using $(a+b)^2\ls 2(a^2+b^2)$) that by the definition (\ref{dist}) we have the following form of the triangle inequality. For $s,t,x\in T$
\be\label{trj}
\varphi_j(s,t)\ls 2(\varphi_j(s,x)+\varphi_j(x,t))
\ee
and as a consequence
\be\label{disjoint}
\forall s,t\in T,\;\varphi_j(s,t)>4a>0\implies B_j(s,a)\cap B_j(t,a)=\emptyset.
\ee
First, let us refine the definition of $j_{n}^{\mu}(t)$ so that they form a non-decreasing sequence which can increase by at most 1. Define
$$\tilde{j}^{\mu}_{n}(t)=\min_{0\ls p\ls n}(j_{p}^{\mu}(t)+n-p).$$
In this way we have $\tilde{j}^{\mu}_{0}(t)=j_{0}^{\mu}(t)$ and  since $j_{p}^{\mu}(t)$ is nondecreasing in $p$ we have for $n\gs 1$ that
\be\label{wzrost}
\tilde{j}^{\mu}_{n}(t)\ls\tilde{j}^{\mu}_{n+1}(t)\ls\tilde{j}^{\mu}_{n}(t)+1.
\ee
Moreover, $\tilde{j}^{\mu}_{n}(t)\ls j_{n}^{\mu}(t)$, so
\be\label{kula1}
\mu(B_{\tilde{j}^{\mu}_{n}(t)}(t,2^n))\gs N_n^{-1}.
\ee 
Finally, for $t\in F$, since $r\gs 4$ we have
\be\label{mainlow2}
\sum_{n\gs0}2^n r^{-\tilde{j}^{\mu}_{n}(t)}\ls\sum_{n\gs0}2^n\sum_{0\ls p\ls n}r^{-j_{p}^{\mu}(t)-n+p}=\sum_{p\gs0}2^p r^{-j_{p}^{\mu}(t)}\sum_{n\gs p}\left(\frac{2}{r}\right)^{n-p}\ls2J_{\mu}(t).
\ee
\noindent The next theorem is a version of Inequality $(4)$ for finite subset $F$ of $T$. It is then extended to the whole $T$ in Theorem \ref{partition1}.
\bt\label{partition}
There exists an admissible sequence of partititons $(\ccA_n)_{n\gs0}$ of $F$ and for $A\in\ccA_n$ there exists an integer $j_n(A)$ such that for each $t\in F$
\be\label{mainlow1}
\sum_{n\gs0}2^n r^{-j_n(A_n(t))}\ls L\sum_{n\gs0}2^n r^{-\tilde{j}^{\mu}_{n}(t)},
\ee
\be\label{pocz}
\mbox{for}\;\; n\ls2, A\in\ccA_n ,j_n(A)=j_0,
\ee
where $j_0$ is as in (\ref{j0}). Moreover,
\be\label{par1'}
A\in\ccA_n,\; C\in\ccA_{n-1},\; A\subset C\implies j_{n-1}(C)\ls j_n(A)\ls j_{n-1}(C)+1.
\ee
and
\be\label{par2}
s,t\in A\in\ccA_n\implies\varphi_{j_n(A)}(s,t)\ls 2^{n+2}.
\ee
\et
\noindent The partitioning procedure is the content of the next lemma.
\begin{lema}\label{par5}
Consider $A\subset F$. There exists a partition $\ccA$ of $A$ such that $|\ccA|\ls N_n$ and for each $B\in\ccA$ there exists $x\in B$ such that
\be\label{par3}
s,t\in B\implies\varphi_{\tilde{j}^{\mu}_{n}(x)}(s,t)\ls 2^{n+4}.
\ee
\end{lema}
\begin{dwd}
Consider $U\subset T$ such that $\forall s,t\in U$ and $s\neq t$, $\varphi_{\tilde{j}^{\mu}_{n}(t)}(s,t)>2^{n+2}$. Let's argue that balls $B_{\tilde{j}^{\mu}_{n}(t)}(t,2^n)$ are disjoint for each $t\in U$. Assume that $\tilde{j}^{\mu}_{n}(s)>\tilde{j}^{\mu}_{n}(t)$ (if $\tilde{j}^{\mu}_{n}(s)=\tilde{j}^{\mu}_{n}(t)$ it follows from (\ref{disjoint})). 
Suppose that there is $x\in B_{\tilde{j}^{\mu}_{n}(s)}(s,2^n)\cap B_{\tilde{j}^{\mu}_{n}(t)}(t,2^n)$. Then we have $\varphi_{\tilde{j}^{\mu}_{n}(s)}(s,x)\ls 2^{n}$ and $\varphi_{\tilde{j}^{\mu}_{n}(t)}(t,x)\ls2^n$. Now, $\tilde{j}^{\mu}_{n}(s)>\tilde{j}^{\mu}_n(t)$ implies
$$\varphi_{\tilde{j}^{\mu}_{n}(t)}(s,x)\ls\varphi_{\tilde{j}^{\mu}_{n}(s)}(s,x),$$
which together with the triangle inequality (\ref{trj}) implies that
$$\varphi_{\tilde{j}^{\mu}_{n}(t)}(t,s)\ls2(\varphi_{\tilde{j}^{\mu}_{n}(t)}(t,x)+\varphi_{\tilde{j}^{\mu}_{n}(t)}(s,x))\ls 2(2^n+ \varphi_{\tilde{j}^{\mu}_{n}(s)}(s,x)) \ls2^{n+2},$$
which contradicts $\varphi_{\tilde{j}^{\mu}_{n}(t)}(s,t)>2^{n+2}$ so $B_{\tilde{j}^{\mu}_{n}(s)}(s,2^n)\cap B_{\tilde{j}^{\mu}_{n}(t)}(t,2^n)=\emptyset$.\\
\noindent By (\ref{kula1}) it follows that $|U|\ls N_n$. Take $U$ with maximal cardinality. Then $$A\subset\bigcup_{t\in U}B_{\tilde{j}^{\mu}_{n}(t)}(t, 2^{n+2})$$ and each of these balls satisfy (\ref{par3}). If we list elements of $U$ denoting them by $t_1,\dots, t_M$, $M\ls N_n$, then the partition $\ccA$ consists of sets 
$$D_1=A\cap B_{\tilde{j}^{\mu}_{n}(t_1)}(t_1, 2^{n+2}),$$ $$D_2=(A\backslash D_1)\cap B_{\tilde{j}^{\mu}_{n}(t_2)}(t,2^{n+2}),$$ $$\vdots$$ 
$$D_M=(A\backslash\bigcup_{k\ls M-1}D_k)\cap B_{\tilde{j}^{\mu}_{n}(t_M)}(t_M,2^{n+2}).$$ Since each element of this partition is contained in the ball of radius $2^{n+2}$ (\ref{par3}) follows.
\end{dwd}
\begin{dwd}[Proof of Theorem \ref{partition}] We proceed by the induction.
Set $\ccA_0=\ccA_1=\ccA_2=\{F\}$ and for $n\ls2$, $j_n(F)=j_0$, so (\ref{pocz}) is satisfied. 
Suppose we have constructed $\ccA_n$ with the following property. For $A\in\ccA_n$, $n\gs2$, there exists an integer, which we denote by $j_n(A)$, such that 
\be\label{par4}
t\in A\implies \tilde{j}^{\mu}_{n-2}(t)=j_n(A).
\ee
 By (\ref{par4}) and (\ref{wzrost}) for $t\in A\in\ccA_n$ we have $\tilde{j}^{\mu}_{n-1}(t)\in\{j_n(A), j_n(A)+1\}$. Set
$$A_0=\{t\in A:\tilde{j}^{\mu}_{n-1}(t)=j_n(A)\}\;\;\mbox{and}\;\;A_1=\{t\in A:\tilde{j}^{\mu}_{n-1}(t)=j_n(A)+1\}.$$
Now, we apply Lemma \ref{par5} for $n-1$ rather than $n$ to get partitions of $A_0$ and $A_1$ into at most $N_{n-1}$ elements, so that we partition the set $A$ into $2N_{n-1}\ls N_n$ sets. For the element $B$ of $A_0$ we put $j_{n+1}(B)=j_n(A)$ and for the element $B$ of $A_1$ we put $j_{n+1}(B)=j_n(A)+1$, so (\ref{par1'}) is satisfied. Apply this procedure to each set $A\in\ccA_{n}$ to get partitition $\ccA_{n+1}$ which is obviously nested. Clearly $|\ccA_{n+1}|\ls N_n^2\ls N_{n+1}$. Condition (\ref{par4}) holds for $n+1$ by the construction as well as (\ref{par2}). Since for $n\ls2$, $j_n(A_n(t))=j_0$ and for $n\gs3$, $j_n(A_n(t))=\tilde{j}^{\mu}_{n-2}(t)$ also (\ref{mainlow1}) follows.
\end{dwd}
\noindent The last step is to level up the partition to the whole set $T$ and to formulate the main lower bound. Before the formal statement let us describe the idea for building the partition elements. We will follow the intuition that the partition of a large finite subset of $T$ should not vary too much from the partition of $T$ itself. To formalize this concept we will consider ascending sequence of finite subsets $F_N\subset T$  and either assign the element of $T$ to already existing partition element of $F_N$ or let it define a new partition element. Alongside, we will need to guarantee that integers $j$ defined for each partition of finite subset converge. This will be done by the procedure of choosing consecutively appropriate subsequences of the indices of $F_N$.
\bt\label{partition1} 
Assume that $T$ is countable. Then, there exists an admissible sequence $(\ccA_n)_{n\gs0}$ of partitions of $T$ and for $A\in\ccA_n$ an integer $j_n(A)$ such that the following holds
\be\label{mainlow3}
\forall t\in T,\;\;\sum_{n\gs0}2^{n}r^{-j_n(A_n(t))}\ls L\E\sup_{t\in T}X_t,
\ee
\be\label{par1}
A\in\ccA_n,\; C\in\ccA_{n-1},\; A\subset C\implies j_{n-1}(C)\ls j_n(A)\ls j_{n-1}(C)+1
\ee
and 
\be\label{par6}
s,t\in A\in\ccA_n\implies\varphi_{j_{n}(A)}(s,t)\ls 2^{n+2}.
\ee
\et
\begin{dwd}
Assume first that $T$ is finite. Then, the result follows from Theorem \ref{partition} combined with (\ref{mainlow2}).\\
Now, let $T$ be countable so that $T=\bigcup_{N\gs1}F_N$, where $(F_N)_{N\gs1}$ is an ascending sequence of finite subsets of $T$. We can enumerate elements of $T$ i.e. $T=\{t_1, t_2,\dots \}$ and put $F_N=\{t_1, t_2,\dots, t_N\}$. Our aim is to construct the admissible sequence of partitions $(\ccA_n)_{n\gs 0}$ of $T$ and verify conditions (\ref{mainlow3}), (\ref{par1}) and (\ref{par6}). This will end the proof. Certainly, $\ccA_0=\{T\}$. The approach is based on the analysis of partitions  $(\ccA_{n, N})_{n\gs0}$ of $F_N$ given by Theorem \ref{partition} and use them for defining the limiting partitions of $T$. Recall that $A_{n,N}(t_i)$ denotes the element of $n$-th partition of $F_N$ that contains $t_i$ and $j_n(A_{n,N}(t_i))$  is the associated integer. To simplify the notation we put $ j_n(A_{n,N}(t_i))=j_{n, N}(t_i)$.\\

\noindent For $t_1$ we obtain a sequence of sets $(A_{n,N}(t_1))_{N\gs 1}$ and the sequence of integers $(j_{n,N}(t_1))_{N\gs 1}$.  Note that by (\ref{par1'})
\be\label{warnaj}
j_0\ls j_{n,N}(t_1)\ls j_{0}+n.
\ee
The main step in the construction is to describe the appropriate limiting procedure allowing to define $(A_n(t_1))_{n\gs0}$ and 
$(j_n(t_1))_{n\gs 0}$. It relies highly on the boundedness of $j_{n, N}(t_i)$ since we can expect a stabilization on certain infinite subsequences. Then the fact that $M$ is a subsequence of $M'$ is equivalent to $M\subset M'$. The first task is to define subsequences $N_n(t_1)$ for $n\gs0$.
We know that for $n\ls 2$ and each $N$, $j_{n, N}(t_1)=j_0$ by (\ref{pocz}), so we put $N_0(t_1)=N_1(t_1)=N_2(t_1)=\{1, 2,\dots\}$. In this case, $\ccA_n=\{T\}$.
Next we choose a subsequence $N_3(t_1)$ of $N_2(t_1)$ such that $(j_{3,N}(t_1))_{N\in N_3(t_1)}$ converges to some limit (which is guaranteed by (\ref{warnaj})) and we denote its' limit by $j_3(t_1)$. 
We proceed in this way to obtain nested sequences $(N_n(t_1))_{n\gs 0}$ and finally by the diagonal procedure we can select a sequence $N(t_1)$ such that for each $n$ the sequence $(j_{n,N}(t_1))_{N\in N(t_1)}$ converges to $j_n(t_1)$. \\
We proceed by induction. Consider $t_i$ for $i>1$ and suppose we have dealt with $t_2,\dots, t_{i-1}$, by which we mean that we have constructed descending family of sets $N(t_2), \dots, N(t_{i-1})$ such that for each $n$ and $j\in\{2,\dots,i-1\}$ the sequence $(j_{n,N}(t_j))_{N\in N(t_j)}$ converges to $j_n(t_j)$  We aim to define an inductive procedure for constructing $N_n(t_i)$ and deciding whether $t_i$ belongs to an already existing partition element or starts a new one. Put $N_0(t_i)=N(t_{i-1})$. Note that $N(t_{i-1})$ is already provided. For $n>0$, consider $N_{n-1}(t_i)$ and suppose that there exists a subsequence $M$ of $N_{n-1}(t_i)$ such that
$$
t_i \in \bigcup^{i-1}_{j=1}A_{n,N}(t_j)\;\; \forall \;\; N\in M.
$$
If so, then we select the smallest $j$ for which there exists $N_n(t_i)\subset N_{n-1}(t_i)$ 
with the property that $t_i\in A_{n,N}(t_j)$ for $N\in N_n(t_j)$. In this way we obtain an infinite subsequence $N_n(t_i)\subset N_{n-1}(t_i)$ and put $t_i$ into the partition element $A_n(t_j)$. Similarly as we have argued for $t_1$, we have that $j_{n,N}(t_i)$ converges to $j_n(t_i)$. Secondly, it may happen that $t_i$ belongs to $\bigcup^{i-1}_{j=1}A_{n,N}(t_j)$ only for finitely many $N\in N_{n-1}(t_i)$. In this case we define a new partition element, $A_n(t_i)$. When choosing $N_n(t_i)\subset N_{n-1}(t_i)$ we care only for the stabilization of $j_{n,N}(t_i)$ namely we require that $j_{n,N}(t_i)$ converges to a limit which we denote by $j_n(t_i)$. Once again, (\ref{warnaj}) implies the existence of this limit. Following the above scheme we decide whether $t_i$ starts a new partition element for $\ccA_n$ or not, construct $j_n(t_i)$ and $N_0(t_i)\supset N_1(t_i)\dots$. We complete the procedure by choosing $N(t_i)$ from the family  $N_0(t_i)\supset N_1(t_i)\dots$ by the diagonal method. Note that it does not affect the convergence of $j_{n,N}(t_i)$ to $j_n(t_i)$ for $N\in N(t_i)$.
\\
Now we check that the defined sequence of partitions is admissible. Namely, that the sequence of partitions $\ccA_n=\{A_n(t_i): i\in I\}$, where $I$ is the index set gathering those points in T which start a partition element,  is nested and satisfies $|\ccA_n|\ls 2^{2^n}$. We have $A_n(t_i)=\{t_i\}\cup\{t_j\in T:\;j>i,\; 
t_j\in A_{n,N}(t_i) \;\forall\; N\in N_n(t_j)\}$. The crucial property of the constructed partition is following. Fix $m$ and consider 
$F_m=\{t_1,...,t_m\}$. For any $n$ there exists a constant $K_{n,m}$ large enough such that for $N> K_{n,m}$ and $N\in N(t_m)$ we have $A_{n,N}(t_i)\supset A_n(t_i) \cap F_m$. Hence,
$|\ccA_n \cap F_m|\ls|\ccA_{n, N}|\ls 2^{2^n}$ and since 
$m$ is arbitrary we conclude that $|\ccA_n|\ls 2^{2^n}$. The fact that $A_{n+1}(t_i)\subset A_{n}(t_i)$ follows clearly from the way we defined those sets.\\
Finally, we verify (\ref{mainlow3}), (\ref{par1}) and (\ref{par6}). They are all straightforward consequences of the fact that for $N>K_{n,m}$ and $N\in N_n(t_m)$ we have $j_{n,N}(t_i)=j_n(t_i)$. 
\end{dwd}
\noindent 
\subsection{Proof of Inequality (5)}
\noindent We want to prove that 
$$\inf\left\{\gamma_2(T_1,d_2)+\gamma_1(T_1, d_{\infty})+\E\sup_{t\in T_2}|X|_t: T\subset T_1+T_2,\;0\in T_1 \right\}{\lesssim}\beta(T, (\varphi_j))+\E\sup_{t\in T}X_t.$$
\noindent Consider $T', T'', T''$ as in Theorem \ref{deco}. Put $T_1=T'$ and $T_2=T''+T'''$.   
Using (\ref{gam2}) and (\ref{gam1}) respectively we get 
\be\label{dol1pop}
\gamma_2(T_1,d_2)\ls L\beta(T, (\varphi_j))),\;\;\gamma_1(T_1, d_\infty)\ls L\beta(T, (\varphi_j)).
\ee
Now,  by replacing $T_2$ by $T_2\cap(T-T_1)$, it follows from Gin\'e-Zinn inequality (\ref{levy1}) that
\begin{eqnarray*}
\E\sup_{t\in T_2}|X|_t&\ls&\sup_{t\in T_2}\int_{\Omega}|t(\omega)|\nu(d\omega)+4\E\sup_{t\in T_2}|X_t|\\
&\ls&\sup_{t\in T'''}\int_{\Omega}|t(\omega)|\nu(d\omega)+L\left(\beta(T, (\varphi_j))+\E\sup_{t\in T-T_1}|X_t|\right),
\end{eqnarray*}
where in the second inequality we used that for $t\in T''$ we have (\ref{L1}). By the triangle inequality $\E\sup_{t\in T-T_1}|X_t|\ls\E\sup_{t\in T}|X_t|+\E\sup_{t\in T_1}|X_t|$.
Now, notice that (\ref{Bern}) together with (\ref{dol1pop}) give that $\E\sup_{t\in T_1}|X_t|\ls L\beta(T, (\varphi_j))$. Hence, the last piece we need to deal with is 
$\sup_{t\in T'''}\int_{\Omega}|t(\omega)|\nu(d\omega)$ Recall that for each $t\in T'''$ there is $s\in T$ such that 
$$|t|\ls 5|s|\1_{\{2|s|\gs r^{-j_0}\}}.$$

\begin{lema}\label{pois1}
Assume $0\in T$ and recall that $r\gs 4$ is fixed. For each $t\in T$ we have
\be\label{pois2}
\int_{\Omega}|t|\1_{\{2|t| \gs r^{-j_0(t)}\}}d\nu\ls L\E\sup_{t\in T}X_t.
\ee
\end{lema}
\begin{dwd}
By (\ref{pois}) we have $\int_{\Omega}|t|\1_{\{2|t| \gs r^{-j_0(t)}\}}d\nu=\E\sum_{i\gs1}|t(Z_i)|\1_{\{2|t(Z_i)| \gs r^{-j_0(t)}\}}$. Define $$N_k=\sum_{i\gs1}\1_{\{r^k\ls2|t(Z_i)|<r^{k+1}\}}$$ and observe that $N_k$ is a Poisson random variable with mean $\nu(\{\omega:r^k\ls2|t(\omega)|<r^{k+1}\})$. Recall that $\varphi_{j_0}(t,0)\ls 4$. Again by (\ref{pois}), $$r^{-2j_0}\varphi_{j_0}(t,0)=\E\sum_{i\gs1}\min(|t(Z_i)|^2,r^{-2j_0}).$$ Hence, 
$$\frac{r^{-2j_0}}{4}\E\sum_{k\gs-j_0}N_k\ls\E\min(|t(Z_i)|^2, r^{-2j_0})\ls r^{-2j_0}.$$
It means that for each $k\gs-j_0$, $\E N_k\ls4$ and since for $\lambda\ls1$ it holds that $\lambda\ls e(1-e^{-\lambda})$ we have
$$\E N_k\ls 4e\P(N_k>0).$$
The above leads to 
$$\E\sum_{i\gs1}|t(Z_i)|\1_{\{2|t(Z_i)| \gs r^{-j_0(t)}\}}\ls\E\sum_{k\gs-j_0}\frac{r^{k+1}}{2}N_k\ls 2e\E\sum_{k\gs-j_0}r^{k+1}\1_{\{N_k>0\}}.$$
Observe that 
\begin{align*}
\E\left(\sum_{k\gs-j_0} r^{2(k+1)}\1_{\{N_k>0\}}\right)^{1/2}&\ls \E\left(\sum_{k\gs-j_0}r^{2(k+1)}N_k\right)^{1/2}\ls L\E\left(\sum_{i\gs1}|t(Z_i)|^2\right)^{1/2}\\
&\ls L\sup_{t\in T}\E\left|\sum_{i\gs1}\va_{i}t(Z_i)\right|\ls L\E\sup_{t\in T}\left|\sum_{i\gs1}\va_{i}t(Z_i)\right|,
\end{align*}
where in the third inequality we used Khintchine's inequality.
In particular, $$\E\left(\sum_{k\gs-j_0} r^{2(k+1)}\1_{\{N_k>0\}}\right)^{1/2}<\infty,$$ so there exists a maximal $k$ for which $N_k>0$. Let's denote it by $K_0$ and observe that it is well-defined on the set $\Omega'=\{\omega: \exists\;\; k\gs -j_0\; \mbox{s.t.}\; N_k>0\}$. We can deduce that the sum is controlled by the last term and write the following
\begin{align*}
\E\sum_{k\gs-j_0}r^{k+1}\1_{\{N_k>0\}}\ls\E\sum_{-j_0\ls k\ls K_0}r^{k+1}\1_{\{N_k>0\}} &\ls\frac{r}{r-1}\E r^{K_0+1}\1_{\Omega'}\\
&\ls\E\left(\sum_{k\gs-j_0} r^{2(k+1)}\1_{\{N_k>0\}}\right)^{1/2}.
\end{align*}
The result follows by (\ref{zero}).
\end{dwd}
\noindent Adding (\ref{pois2}) to the discussion at the beginning of this section and applying  (\ref{zero}) yields Inequality $(5)$.

\section{Appendix}

\subsection{Extension to the separable case}\label{A1}
Notice that the proof of Theorem \ref{partition1} can be extended to separable $T$, that is to such $T$ which has a countable subset which is dense in $T$ with respect to some metric and also each function $\varphi_j(s,t)$ is continuous in this metric. Let $\ccA_n$ be the admissisble sequence of partitions of this subset given by Theorem \ref{partition1}. Let $\bar{A}$ be the closure of the set $A$, so by the continuity of $\varphi_{j(A)}(s,t)$, $A\in\ccA_n$, we have for each $s,t\in \bar{A}$ that $\varphi_{j(A)}(s,t)\ls2^{n+2}$. Sets $\bar{A}$ cover $T$, but they do not have to be disjoint. We construct an admissible sequence $\ccB_n$ of $T$ such that for each $B\in\ccB_n$ there exists $A\in\ccA_n$ such that $B\subset\bar{A}$. It is possible because $\ccA_n$ is admissible, so in particular $\bar{A}=\bigcup\{\bar{C}:\; C\subset A,\; A\in\ccA_{n+1}\}$.

\subsection{Symmetrization}\label{A2}
Suppose we are given some non-Gaussian infinitely divisible process $(X_t)$ and we want to apply Theorem \ref{mainpop}. Let $(X_t')$ be an independent copy of the process $(X_t)$. First of all, it can be applied only if
$$A:=\E\sup_{t\in T}(X_t-X_t')<\infty.$$
Then,
$$\E\sup_{t\in T}|X_t-X_t'-X_s+X_s'|=2A,$$
so by Fubini's theorem
$$\sup_{s,t\in T}\E|X_t-X_s|<\infty,$$
which proves that $\sup_{t\in T}\E|X_t|<\infty$ (since we can assume $0\in T$) and 
$$2A\gs\E\sup_{t\in T}|X_t-\E X_t|.$$
On the other hand,
$$A\ls\E\sup_{t\in T}|X_t-\E X_t|+\E\sup_{t\in T}|X_t'-\E X_t'|=2\E\sup_{t\in T}|X_t-\E X_t|,$$
so $A\asymp \E\sup_{t\in T}|X_t-\E X_t|=:B$.
Finally,
$$\E\sup_{t\in T}|X_t|\ls B+ \sup_{t\in T}|\E X_t|\ls 2\max\{B, \sup_{t\in T}|\E X_t|\}.$$
Now, to show the reverse inequality notice that
$$\E\sup_{t\in T}|X_t|\gs\sup_{t\in T}\E |X_t|\gs\sup_{t\in T}|\E X_t|$$
and thus
$$\E\sup_{t\in T}|X_t|+\sup_{t\in T}|\E X_t|\gs\E\sup_{t\in T}(|X_t|+|\E X_t|)\gs\E\sup_{t\in T}|X_t-\E X_t|.$$
It shows that $2\E\sup_{t\in T}|X_t|\gs B$ and finally
$$2\E\sup_{t\in T}|X_t|\gs\max\{B, \sup_{t\in T}|\E X_t|\}.$$

\subsection{Decomposition of the process}\label{A3}
We can deduce the following version of Theorem \ref{mainpop}
\begin{theo}
Consider an infinitely divisible process $(X_t)_{t\in T}$ and assume $T$ is countable set. Then 
$$\inf\{\gamma_2(T,d_2^1)+\gamma_1(T,d_{\infty}^1)+\E\sup_{t\in T}|X^1|_t: X_t=X_t^1+X_t^2\}\asymp\E\sup_{t\in T}X_t,$$
where $(X_t^1)_{t\in T}$, $(X_t^2)_{t\in T}$ are infinitely divisible and distances $d_2^1$ and $d_{\infty}^1$ are as in (\ref{charac}) but for the measure $\nu^1$ associated with the process $(X_t^1)_{t\in T}$. 
\end{theo}

\begin{dwd}
The upper bound is obvious since $\E\sup_{t\in T}X_t\ls\E\sup_{t\in T}X_t^1+\E\sup_{t\in T}X_t^2$. For the first term we apply chaining bound (\ref{Bern}) while for the second the bound (\ref{Gora}).\\
\noindent For the lower bound we simply use the decomposition in Theorem \ref{mainpop}. Since every point $t\in T$ can be written as $t=t_1+t_2$, where $t_1\in T_1$ and $t_2\in T_2$, we set $X_t^1=X_{t_1}$ and $X_t^2=X_{t_2}$, which finishes the proof.
\end{dwd}

\section*{Acknowledgements}
The authors would like to thank prof. Michel Talagrand for the most helpful comments and suggestions on the whole approach to the problem, especially for pointing out that proofs of steps $(5)$ and $(4)$ are possible in their present form.\\
\noindent Also, we are very grateful to the anonymous referees for suggesting substantial changes of the exposition of the material in the paper which made it significantly clearer and more readable as well as for pointing out numerous errors in the initial version of the paper.\\ 
\noindent The first author was supported by NCN Grant UMO-2016/21/B/ST1/01489.\\
\noindent The second author was supported by NCN Grant UMO-2018/31/N/ST1/03982.

%%%%%%%%%%%%%%%%%%%%%%%%%%%%%%%%%%%%%%%%%%%%%%
%% Supplementary Material, if any, should   %%
%% be provided in {supplement} environment  %%
%% with title and short description.        %%
%%%%%%%%%%%%%%%%%%%%%%%%%%%%%%%%%%%%%%%%%%%%%%

%%  MR numbers will be added by VTeX.                      %%
%%                                                         %%
%%  Use \cite{...} to cite references in text.             %%
%%                                                         %%
%%%%%%%%%%%%%%%%%%%%%%%%%%%%%%%%%%%%%%%%%%%%%%%%%%%%%%%%%%%%%

%% if your bibliography is in bibtex format, uncomment commands:
%\bibliographystyle{imsart-number} % Style BST file (imsart-number.bst or imsart-nameyear.bst)
%\bibliography{bibliography}       % Bibliography file (usually '*.bib')

\begin{thebibliography}{4}
\bibitem{Bed2} Bednorz, W.: The majorizing measure approach to sample boundedness, \emph{Colloquium Mathematicum} {\bf 139} (2015), 205--227.

\bibitem{Bed1} Bednorz, W. and Lata\l{}a, R.: On the boundedness of Bernoulli processes,  \emph{Ann. of Math.} {\bf 180} 
(2014), 1167--1203.

\bibitem{Dal}Daley, D.J. and Vere-Jones, D.: \textit{An introduction to the theory of point processes. Vol. 1.} (2003), Springer-Verlag, New York.

\bibitem{Fer} Fernique, X.: R{\'e}gularit{\'e} des trajectoires des fonctions al{\'e}atoires gaussiennes,
in: \textit{{\'E}cole d'{\'E}t{\'e} de Probabilit{\'e}s de Saint-Flour, IV-1974,} \emph{Lecture Notes in Math.} {\bf 480},
 Springer, Berlin, 1975, 1--96.

\bibitem{von} Kneser, H.: \textit{Sur un th\'eor\'eme fondamental de la th\'eorie des jeux}. In: Hofmann, K. and Betsch, G. (eds.) Gesammelte Abhandlungen / Collected Papers. Berlin, New York De Gruyter, (2011),  484--486.

\bibitem{Last} Last, G. and  Penrose, M.: \textit{Lectures on the Poisson Process}. Institute of Mathematical Statistics Textbooks. Cambridge Univ. Press, Cambridge, 2018.



\bibitem{M-R} Marcus, M. B. and Rosen, J.: Markov processes, Gaussian processes, and local times,  \emph{Cambridge Stud. Adv. Math.} {\bf 100},  Cambridge Univ. Press, Cambridge, 2006.

\bibitem{Ros3} Marcus, M.B., Rosi\'{n}ski, J.: \textit{Sufficient conditions for boundedness of moving averages.}In: Stochastic Inequalities and Applications. Progr. Probab., vol. {\bf 56}, Birkh\"{a}user, Basel (2003), 113--128.

\bibitem{Ros1}Rosi\'{n}ski, J.: \textit{On series representations of infinitely divisible random vectors.} Ann. Probab., {\bf18}, (1990), pp. 405--430.

\bibitem{Ros4}Rosi\'{n}ski, J.: \textit{Series representation of L\'evy processes from the perspective of point processes.} In: Barndnoff-Nielsen, O.E., Mikosch, T., Resnick, S.I. (eds.) L\'evy Processes Theory and Applications, Birkh\"{a}user, Basel (2001), 401--415.

\bibitem{Ros2}Rosi\'{n}ski, J.: \textit{Representations and isomorphism identities for infinitely divisible processes.} Ann. Probab., {\bf46}, (2018), pp. 3229--3274.

\bibitem{Sat1}Sato, K.:\textit{L\'{e}vy Processes and Infinitely Divisible Distributions.} Cambridge Studies in Advanced Mathematics, {\bf68}, Cambridge University Press, (2013), Cambridge.

\bibitem{Tal} Talagrand, M. (2019)  Upper and Lower Bounds for Stochastic Processes. Modern Methods and Classical Problems. Preprint.

\bibitem{Tal1} Talagrand, M.:  Upper and Lower Bounds for Stochastic Processes. Modern Methods and Classical Problems,
  \emph{Ergeb. Math. Grenzgeb.} {\bf 60}. (2014) Springer, New York. 
  
\bibitem{Tal5} Talagrand, M.: Regularity of Gaussian processes. \emph{Acta Math.}, {\bf 159} (1987), 99--149.

\bibitem{Tal2} Talagrand, M.: Regularity of infinitely divisible processes. \emph{Ann. Probab.} {\bf 20}. Number 1 (1993), 362--432. 

\bibitem{Tal3} Talagrand, M.: Majorising measures without measures. \emph{Ann. Probab.} {\bf 29}. Number 1 (2001), 411--417.

\bibitem{Tal4} Talagrand, M.: Sample Boundedness of Stochastic Processes Under Increment Conditions. \emph{Ann. Probab.} {\bf 18}. Number 1 (1990), 1--49.


\end{thebibliography}

%% or include bibliography directly:

\end{document}